\documentclass[a4paper,10pt]{article}
\usepackage{a4}

\newtheorem{definition}{\hspace {1cm}Definition}[section]
\newtheorem{proposition}{\hspace {1cm}Proposition}[section]
\newtheorem{theorem}{\hspace {1cm}Theorem}[section]
\newtheorem{remark}{\hspace {1cm}Remark}[section]
\newtheorem{example}{\hspace {1cm}Example}[section]

\begin{document}

\author{Adelina Manea}
\title{A de Rham theorem with respect to the Liouville foliation on $TM^0$, for a Finsler manifold $M$}
\date{}
\maketitle

\begin{abstract}
On the slit tangent manifold of a Finsler manifold there are given the vertical and the Liouville foliations, \cite{Bej}. In this paper we define new types of vertical forms with respect to the Liouville foliation. We prove a de Rham type theorem using this forms.  
\end{abstract}


\textbf{AMS 2000}: 53C12, 53C60.

\mathstrut

\section{Preliminaries}

\setcounter{equation}{0}  Finding new topological invariants of differentiable manifolds is still an open problem for geometries. The cohomology groups are such invariants. The Finsler manifolds are interesting models for some physical phenomena, so their properties are also useful to investigate, \cite{Bej}, \cite{Mir}. The cohomology groups of manifolds, related sometimes to some foliations on them, have been studied in the last decades, \cite{T}, \cite{V}, \cite{Bull}. Our present work intends to develop the study of the Finsler manifolds and the foliated structures of the tangent bundle of such a manifold.

For the beginning, we present two foliations on the slit tangent manifold $TM^0$ of a $n$-dimensional Finsler manifold $(M,F)$, following \cite{Bej}. In this paper the indices take the values $i,j,i_1,j_1,...$ $=\overline{1,n}$ and $a,b,a_1,b_1,...$$=\overline{1,n-1}$.

Let $(M,F)$ be a $n$-dimensional Finsler manifold and $G$ the Sasaki-Finsler metric on its slit tangent manifold $TM^0$. The vertical bundle $VTM^0$ of $TM^0$ is the tangent (structural) bundle to the vertical foliation $F_V$ determined by fibers $\pi :TM^0\rightarrow M$. If $(x^i,y^i)_{i=\overline{1,n}}$ are local coordinates on $TM^0$, then $VTM^0$ is locally spanned by $\{\frac{\partial}{\partial y^i}\}_i$. A canonical transversal (also called horizontal) distribution is constructed in \cite{Bej} as follows. We denote by $(g^{ij}(x,y))_{i,j}$ the inverse matrix of $g=(g_{ij}(x,y))_{i,j}$, where 

\begin{eqnarray}
	g_{ij}(x,y)=\frac{1}{2}\frac{\partial^2 F^2}{\partial y^i \partial y^j}(x,y),
\end{eqnarray}
and $F$ is the fundamental function of the Finsler manifold. Obviously, we have the equalities $\frac{\partial g_{ij}}{\partial y^k}$$=\frac{\partial g_{ik}}{\partial y^j}$$=\frac{\partial g_{jk}}{\partial y^i}$. 

Then locally define the functions
	\[G^i=\frac{1}{4}g^{ik}\left(\frac{\partial^2 F^2}{\partial y^k \partial x^h}y^h-\frac{\partial F^2}{\partial x^k}\right),\quad G^j_i=\frac{\partial G^j}{\partial y^i}.
\]
There exists on $TM^0$ a $n$ distribution $HTM^0$ locally spanned by the vector fields
\begin{eqnarray}
	\frac{\delta }{\delta x^i}=\frac{\partial }{\partial x^i}-G_i^j\frac{\partial }{\partial y^j},\quad (\forall)i=\overline{1,n}.
\end{eqnarray}
The Riemannian metric $G$ on $TM^0$ is satisfying 
\begin{eqnarray}
	G(\frac{\delta}{\delta x^i},\frac{\delta}{\delta x^j})=G(\frac{\partial}{\partial y^i},\frac{\partial}{\partial y^j})=g_{ij},\quad G(\frac{\delta}{\delta x^i},\frac{\partial}{\partial y^j})=0,\quad (\forall)i,j.
\end{eqnarray}
The local basis $\{\frac{\delta}{\delta x^i},\frac{\partial}{\partial y^i}\}_i$ is called adapted to the vertical foliation $F_V$ and we have the decomposition
\begin{eqnarray}
	TTM^0=HTM^0\oplus VTM^0.
\end{eqnarray}

Now, let $Z$ be the global defined vertical Liouville vector field on $TM^0$,
\begin{eqnarray}
	Z=y^i\frac{\partial}{\partial y^i},
\end{eqnarray}
and $L$ the space of line fields spanned by $Z$. We call this space \textit{the Liouville distribution} on $TM^0$. The complementary orthogonal distributions to $L$ in $VTM^0$ and $TTM^0$ are denoted by $L'$ and $L^\perp$, respectively. It is proved, \cite{Bej}, that the both distributions $L'$ and $L^\perp$ are integrable and we also have the decomposition
\begin{eqnarray}
	VTM^0=L'\oplus L.
\end{eqnarray}
Moreover, we have, \cite{Bej}:
\begin{proposition}
a) The foliation determined by the distribution $L^\perp$ is just the foliation determined by the level hypersurfaces of the fundamental function $F$ of the Finsler manifold.

b) For every fixed point $x_0\in M$, the leaves of the Liouville foliation $F_{L'}$ determined by the distribution $L'$ on $T_{x_0}M$ are just the $c$-indicatrices of $(M,F)$:
\begin{eqnarray}
	I_{x_0}M(c):\quad F(x_0,y)=c,\quad (\forall)y\in T_{x_0}M.
\end{eqnarray}

c) The foliation $F_{L'}$ is a subfoliation of the vertical foliation.
\end{proposition}

As we already saw, the vertical bundle is locally spanned by $\{\frac{\partial}{\partial y^i}\}_{i=\overline{1,n}}$ and it admits decomposition (1.6). In the following we give another basis on $VTM^0$, adapted to $F_{L'}$. 

There are some useful facts which follow from the homogeneity of the fundamental function of the Finsler manifold $(M,F)$. By the Euler theorem on positively homogeneous functions we have, \cite{Bej},
\begin{eqnarray}
	F^2(x,y)=y^iy^jg_{ij}(x,y),\quad \frac{\partial F}{\partial y^k}=\frac{1}{F}y^ig_{ki}, \quad y^i\frac{\partial g_{ij}}{\partial y^k}=0,\quad \forall k=\overline{1,n}.
\end{eqnarray}
Hence it results
\begin{eqnarray}
	G(Z,Z)=F^2.
\end{eqnarray}

We consider the following vertical vector fields:
\begin{eqnarray}
	X_k=\frac{\partial}{\partial y^k}-t_kZ,\quad k=\overline{1,n},
\end{eqnarray}
where the functions $t_i$ are defined by the conditions
\begin{eqnarray}
	G(X_k,Z)=0, \forall k=\overline{1,n}.
\end{eqnarray}
The above conditions become
	\[G(\frac{\partial}{\partial y^k}, y^i\frac{\partial}{\partial y^i})-t_kG(Z,Z)=0,
\]
so, taking into account also (1.3) and (1.9), we obtain the local expression of functions $t_k$ in a local chart $(U,(x^i,y^i))$:
\begin{eqnarray}
	t_k=\frac{1}{F^2}y^ig_{ki}=\frac{1}{F}\frac{\partial F}{\partial y^k},\quad \forall k=\overline{1,n}.
\end{eqnarray}

If $(\tilde{U},(\tilde{x}^{i_1},\tilde{y}^{i_1}))$ is another local chart on $TM^0$, in $U \cap \tilde{U} \neq \oslash $ we have:
	\[\tilde{t}_{k_1}=\frac{1}{F^2}\tilde{y}^{i_1}\tilde{g}_{i_1k_1}=\frac{1}{F^2}\frac{\partial \tilde{x}^{i_1}}{\partial x^i}y^i \frac{\partial x^k}{\partial \tilde{x}^{k_1}}\frac{\partial x^i}{\partial \tilde{x}^{i_1}}g_{ki}=\frac{\partial x^k}{\partial \tilde{x}^{k_1}}t_k,
\]
so we obtain the following changing rule for the vector fields (1.10):
\begin{eqnarray}
	\tilde{X}_{i_1}=\frac{\partial x^k}{\partial \tilde{x}^{i_1}}X_k,\quad \forall i_1=\overline{1,n}.
\end{eqnarray}
By a straightforward computation, using (1.8), it results:
\begin{proposition}
\quad The functions $\{t_k\}_{k=\overline{1,n}}$ defined by (1.12) are satisfying:
\begin{eqnarray}
\it{a})	\quad y^it_i=1;\quad y^i X_i=0;
\end{eqnarray}
\begin{eqnarray}
	\it{b})\quad \frac{\partial t_l}{\partial y^k}=-2t_kt_l+\frac{1}{F^2} g_{kl},\quad Zt_k=-t_k,\quad \forall k,l=\overline{1,n};
\end{eqnarray}
\begin{eqnarray}
	d)\quad y^j\frac{\partial t_j}{\partial y^i}=-t_i,\quad \forall i=\overline{1,n},\quad y^i(Zt_i)=-1,\quad y^i(ZX_i)=0.;
\end{eqnarray}
\end{proposition}
\begin{proposition}
There are the relations:
\begin{eqnarray}
	[X_i,X_j]=t_iX_j-t_jX_i,
\end{eqnarray}
\begin{eqnarray}
	[X_i,Z]=X_i,
\end{eqnarray}
for all $i,j=\overline{1,n}$.
\end{proposition}
By the conditions (1.11), $\{X_1,...,X_n\}$ are $n$ vector fields orthogonal to $Z$, so they belong to the $(n-1)$-dimensional distribution $L'$. It results that they are linear dependent and, from (1.14), 
\begin{eqnarray}
	X_n=-\frac{1}{y^n}y^aX_a,
\end{eqnarray}
since the local coordinate $y^n$ is nonzero everywhere.

We also proved that, \cite{Manea}:

\begin{proposition}
The system of vector fields $\{X_1,X_2,...,X_{n-1},Z\}$ of vertical vector fields is a locally adapted basis to the Liouville foliation $F_{L'}$, on $VTM^0$.
\end{proposition}

The whole proofs for propositions from this section are given in \cite{Manea}.

More clearly, let $(\tilde{U},(\tilde{x}^{i_1},\tilde{y}^{i_1}))$, $(U,(x^i,y^i))$ be two local charts which domains overlap, where $\tilde{y}^k$ and $y^n$ are nonzero functions (in every local charts on $TM^0$ there is at least one nonzero coordinate function $y^i$). The adapted basis in $\tilde{U}$ is $\{\tilde{X}_1,\tilde{X}_2,...,\tilde{X}_{k-1}, \tilde{X}_{k+1},...,\tilde{X}_n,Z\}$. In $U \cap \tilde{U}$ we have (1.13) and (1.19), hence
\[\tilde {X}_{i_1}=\sum_{i=1}^{n-1}(\frac{\partial x^i}{\partial \tilde{x}^{i_1}}-\frac{y^i}{y^n}\frac{\partial x^n}{\partial \tilde{x}^{i_1}})X_i;\quad X_j=\sum_{j_1=1, j_1 \neq k}^n(\frac{\partial \tilde{x}^{j_1}}{\partial x^j}-\frac{\tilde{y}^{j_1}}{\tilde{y}^k}\frac{\partial \tilde{y}^k}{\partial x^j})\tilde{X}_{j_1},
\]
for all $i_1=\overline{1,n}$, $i_1\neq k$, $j=\overline{1,n-1}$. Ones can see that the above relations also imply
\[ \frac{\partial x^i}{\partial \tilde{x}^{k}}-\frac{y^i}{y^n}\frac{\partial x^n}{\partial \tilde{x}^{k}}=-\sum_{i_1=1,i_1 \neq k}^{n}\frac{\tilde{y}^{i_1}}{\tilde{y}^k}(\frac{\partial x^i}{\partial \tilde{x}^{i_1}}-\frac{y^i}{y^n}\frac{\partial x^n}{\partial \tilde{x}^{i_1}}).
\]
By a straightforward calculation we have that the determinant of the changing matrix $\{X_1,X_2,...,X_{n-1},Z\}$$ \rightarrow$ $\{\tilde{X}_1,\tilde{X}_2,...,\tilde{X}_{k-1}, \tilde{X}_{k+1},...,\tilde{X}_n,Z\}$ on $L'$ is equal to 
\[ (-1)^{n+k}\frac{\tilde{y}^k}{y^n} det\left(\frac{\partial x^i}{\partial \tilde{x}^j}\right)_{i,j=\overline{1,n}}.
\]

\section{New types of vertical forms with respect to Liouville foliation on $TM^0$}
\setcounter{equation}{0}
Now, let $\{\delta y^i=dy^i+G_j^idx^j\}_{i=\overline{1,n}}$ be the dual basis of $\{\frac{\partial}{\partial y^i}\}_{i=\overline{1,n}}$ on $VTM^0$. We also consider the space $\Omega^0(TM^0)$ of differentiable functions on $TM^0$, the module $\Omega^{0,q}(TM^0)$ of vertical $q$-forms and the foliated derivative $d_{01}$ with respect to the vertical foliation on $TM^0$.
\begin{proposition}
The vertical $1$-form $\omega_0=t_i\delta y^i$ is globally defined and 
\begin{eqnarray}
	\omega_0(Z)=1,\quad \omega_0(X_a)=0,\quad \omega_0=d_{01}(lnF),
\end{eqnarray}
for all $a=\overline{1,n-1}$, $X_a$ given by (1.10) and $F$ the fundamental function of the Finsler manifold.
\end{proposition}

\textit{Proof:} In $\tilde{U}\cap U$ we have
\[\tilde{\omega}_0=\tilde{t}_{i_1}\delta \tilde{y}^{i_1}=\frac{\partial x^i}{\partial \tilde{x}^{i_1}}t_i\frac{\partial \tilde{x}^{i_1}}{\partial x^j}\delta y^j =t_i\delta y^i =\omega_0.
\]
We also have $\delta y^i(Z)=y^i$, for all $i=\overline{1,n}$, and taking into account the first relation (1.14), \[\omega_0(Z)=1,\quad \omega_0(X_a)=t_i \delta y^i(\frac{\partial }{\partial y^a}-t_aZ)=t_i\delta^i_a-t_at_iy^i=0,\]
where $\delta^i_a$ is the Kronecker symbol. By the relation (1.19) it results also $\omega_0(X_n)=0$. In the last, locally we have
	\[d_{01}(lnF)=\frac{\partial (lnF)}{\partial y^i}\delta y^i=\frac{1}{F}\frac{\partial F}{\partial y^i}\delta y^i=\omega_0,
\]where we used relation (1.12).
The equality $\omega_0=d_{01}(lnF)$ shows that $\omega_0$ is a $d_{01}$-exact vertical 1-form and the Liouville foliation $L'$ is defined by the equation $\omega_0=0$.
\begin{definition}
A vertical $q$-form $\omega \in \Omega ^{0,q}(TM^0)$ is called a vertical $(s,t)$-form or a $(0,s,t)$-form iff for vertical vector fields $Y_1,Y_2,...,Y_q$, $\omega(Y_1,...,Y_q)\neq 0$ only if $s$ arguments are in $L'$ and $t$ arguments are in $L$.
\end{definition}
Since $L$ is a line distribution, we can talk only about $(0,s,t)$-forms with $t\in \{0,1\}$. We denote the space of $(0,s,t)$-forms by $\Omega^{0,s,t}(TM^0)$. By the above definition, we have the equivalence
\begin{eqnarray}
	\omega \in \Omega^{0,q-1,1}(TM^0)\quad \Longleftrightarrow \quad \omega(Y_1,...,Y_q)= 0,\quad (\forall)Y_1,...,Y_q \in \{X_1,...,X_{n-1}\},
\end{eqnarray}
where $\{X_i\}_{i=\overline{1,n-1}}$ is the local basis in $L'$ from proposition 1.4.
\begin{proposition}
Let $\omega$ be a nonzero vertical $q$-form. The following assertions are true:

a) $\omega \in \Omega ^{0,q,0}(TM^0)$ iff $i_Z\omega =0$, where $i_Z$ is the interior product with the vertical Liouville vector field $Z$.

b) The vertical $(q-1)$-form $i_Z\omega$ is a $(0,q-1,0)$-form.

c) $\omega \in \Omega ^{0,q-1,1}(TM^0)$ implies $i_Z\omega \neq 0$.

d) If there is a $(0,q-1,0)$-form $\alpha$ such that $\omega=\omega_0\wedge \alpha$, then $\omega \in \Omega ^{0,q-1,1}(TM^0)$.
\end{proposition}
\textit{Proof:}a) Let  $\omega \in \Omega ^{0,q,0}(TM^0)$, hence $\omega(Y_1,...,Y_q)\neq 0$ only if all the arguments are in $L'$. So, $i_Z\omega $ is a $(q-1)$-form and $i_Z\omega(Y_1,...,Y_{q-1})=\omega (Z,Y_1,...,Y_{q-1})=0$, for every vertical vector fields $Y_1,...,Y_{q-1}$. That means $i_Z\omega=0$. Conversely, if $\omega$ is a $(0,q)$-form such that $i_Z\omega=0$, then $\omega(Y_1,...,Y_q)= 0$ since there is an index $i\in \{1,..., q\}$ such that $Y_i=Z$. Hence $\omega$ does not vanish only on $L'$, and by definition it is a $(0,q,0)$-form.

b) We have $i_Zi_Z\omega=0$, from the definition of a form on a manifold, and taking into account a), it results that $i_Z\omega$ is a $(0,q-1,0)$-form.

c) If $\omega$ is a nonzero $(0,q-1,1)$-form, then $\omega(Y_1,...,Y_q)\neq 0$ only if exactly one of the arguments is from the line distribution $L=span{Z}$. Then $i_Z\omega ((Y_1,...,Y_{q-1})\neq 0$ for some vertical vector fields $Y_1,...,Y_{q-1}\in L'$.

d) Let $\alpha$ be a form like in hypothesis, and $Y_1,...,Y_q$, $q$ arbitrary vertical vector fields. 
	\[ \omega(Y_1,...,Y_q)=(\omega_0 \wedge \alpha)(Y_1,...,Y_q)=\sum_{\sigma \in S_q}\epsilon(\sigma) \omega_0(Y_{\sigma (1)})\alpha (Y_{\sigma (2)},...,Y_{\sigma (q)}).
\]
But $\omega_0$ vanishes on $L'$ and for the all arguments $Y_1,...,Y_q$ in $L'$, the above sum has all terms nulls. Taking into account relation (2.2), we have $\omega \in \Omega^{0,q-1,1}(TM^0)$.

\begin{proposition}
For every vertical $q$-form $\omega$ there are $\omega_1\in \Omega^{0,q,0}(TM^0)$ and $\omega_2\in \Omega^{0,q-1,1}(TM^0)$ such that $\omega=\omega_1+\omega_2$, uniquely.
\end{proposition}
\textit{Proof:} Let $\omega$ be a nonzero vertical $q$-form. If $i_Z\omega =0$, then $\omega \in \Omega^{0,q,0}(TM^0)$ from Proposition 2.2, so $\omega =\omega +0$.
If $i_Z\omega \neq 0$, then let $\omega_2$ be the vertical $q$-form $\omega_0\wedge i_Z\omega$. By Proposition 2.2 d), it results $\omega_2$ is a $(0,q-1,1)$-form. Moreover, putting $\omega_1=\omega-\omega_2$, 
\begin{eqnarray}
	i_Z\omega_1=i_Z\omega-i_Z(\omega_0\wedge i_Z\omega)=i_Z\omega-\omega_0(Z)i_Z\omega=0,
\end{eqnarray}
where we used relation(2.1). So, $\omega_1$ is a $(0,q,0)$-form and $\omega_1$, $\omega_2$ are unique defined by $\omega$. Obviously $\omega=\omega_1+\omega_2$. 

 We have to remark that only the zero $q$-form could be a $(0,q,0)$- and a $(0,q-1,1)$-form at the same time. The proposition 2.3 proves the decomposition
 
\begin{eqnarray}
	\Omega^{0,q}(TM^0) = \Omega^{0,q,0}{TM^0} \oplus  \Omega^{0,q-1,1}(TM^0).
\end{eqnarray}

A consequence of the propositions 2.2 and 2.3 is:

\begin{proposition}
Let $\omega$ be a $(0,q)$-form. We have the equivalence:
\begin{eqnarray}
	\omega \in \Omega^{0,q-1,1}(TM^0)\quad \Longleftrightarrow \quad (\exists)\alpha \in \Omega^{0,q-1,0}(TM^0),    
	\omega=\omega_0\wedge \alpha.
\end{eqnarray}
\end{proposition}

Taking into account the characterization given in Proposition 2.2a) and relation (2.5), ones can see that:

\begin{proposition}
We have the following facts:

a) If $\omega\in \Omega^{0,q,0}(TM^0)$ and $\theta \in \Omega^{0,r,0}(TM^0)$, then $\omega\wedge \theta \in \Omega^{0,q+r,0}(TM^0)$.

b) If $\omega\in \Omega^{0,q,1}(TM^0)$ and $\theta \in \Omega^{0,r,0}(TM^0)$, then $\omega\wedge \theta \in \Omega^{0,q+r,1}(TM^0)$.

c)If $\omega\in \Omega^{0,q,1}(TM^0)$ and $\theta \in \Omega^{0,r,1}(TM^0)$, then $\omega\wedge \theta =0$.
\end{proposition}

\begin{example}
a) $\omega_0$ is a $(0,0,1)$-form because there is the $(0,0,0)$-form, the constant 1 function on $TM^0$, such that $\omega_0=\omega_0\cdot 1$.

b) $\theta_i=\delta y^i-y^i\omega_0$ is a $(0,1,0)$-form, for each $i-\overline{1,n}$. Indeed, 
	\[\theta_i(Z)=\delta y^i(Z)-\omega_0(Z)y^i=0,
\]so $i_Z\theta_i=0$. We have to remark that the vertical $1$-forms $\{\theta_i\}_{i=\overline{1,n}}$ are linear dependent, since $\sum t_i\theta_i=0$.

c) $i_Z(\theta_i\wedge \theta_j)(Y)=\theta_i(Z)\theta_j(Y)-\theta_j(Z)\theta_i(Y)=0$, for every vertical vector field $Y$, hence $\theta_i\wedge \theta_j \in \Omega^{0,2,0}(TM^0)$.
\end{example}

\begin{proposition}
The foliated derivative $d_{01}:\Omega^{0,q}(TM^0)\rightarrow\Omega^{0,q+1}(TM^0)$ has the following property: for every $(0,q-1,1)$-form $\omega$, $d_{01}\omega$ is a $(0,q,1)$-form. 
\end{proposition}

\textit{Proof:} Let $\omega$ be a $(0,q-1,1)$-form. From relation (2.5), there is a $(0,q-1,0)$ form $\alpha$ such that $\omega=\omega_0\wedge \alpha$. Proposition 2.3 proved also that $\alpha=i_Z\omega$. Taking into account that $\omega_0$ is an $d_{01}$-exact form, it follows
	\[d_{01}\omega=d_{01}(\omega_0\wedge \alpha)=-\omega_0\wedge d_{01}\alpha=-\omega_0\wedge \beta_1-\omega_0\wedge\beta_2,
\]with $\beta_1$, $\beta_2$ the $(0,q,0)$-, $(0,q-1,1)$-forms, components of the $(0,q)$-form $d_{01}\alpha$. But $\beta_2=\omega_0\wedge \theta$ from relation (2.5), so we have $d_{01}\omega=-\omega_0\wedge \beta_1$; it follows $d_{01}\omega \in \Omega^{0,q,1}(TM^0)$. We can write 
\begin{eqnarray}
	d_{01}(\Omega^{0,q-1,1}(TM^0)) \subset \Omega^{0,q,1}(TM^0).
\end{eqnarray}
Let as consider $\xi_1$, $\xi_2$ the projections of the module $\Omega^{0,q}(TM^0)$ on its direct summands from relation (2.4). 
\begin{eqnarray}
	\xi_1:\Omega^{0,q}(TM^0)\rightarrow \Omega^{0,q,0}(TM^0),\quad \xi_1(\omega)=\omega-\omega_0\wedge i_Z\omega,\quad (\forall)\omega\in \Omega^{0,q}(TM^0),
\end{eqnarray}
\begin{eqnarray}
	\xi_2:\Omega^{0,q}(TM^0)\rightarrow \Omega^{0,q-1,1}(TM^0),\quad \xi_2(\omega)=\omega_0\wedge i_Z\omega,\quad (\forall)\omega\in \Omega^{0,q}(TM^0),
\end{eqnarray}

\begin{remark}
For an arbitrary $(0,q)$-form $\omega$, $d_{01}\omega=d_{01}(\xi_1(\omega))+d_{01}((\xi_2\omega))$. Relation (2.6) show that $d_{01}(\xi_2(\omega))$ is a $(0,q,1)$-form, hence $\xi_1(d_{01}(\xi_2(\omega)))=0$. It results 
\begin{eqnarray}
	\xi_1(d_{01}\omega)=\xi_1(d_{01}(\xi_1(\omega),\quad \xi_2(d_{01}\omega)=\xi_2(d_{01}(\xi_1(\omega)+d_{01}(\xi_2(\omega)).
\end{eqnarray}
The above relations prove that $d_{01}(\Omega^{0,q,0}(TM^0)) \subset \Omega^{0,q+1,0}(TM^0)\oplus \Omega^{0,q,1}(TM^0)$.
\end{remark}

Let us define the following operators:
\begin{eqnarray}
	d':\Omega^{0,q,0}(TM^0)) \rightarrow \Omega^{0,q+1,0}(TM^0),\quad d'(\omega)=\xi_1(d_{01}\omega),
\end{eqnarray}

\begin{eqnarray}
	d":\Omega^{0,q,0}(TM^0)) \rightarrow \Omega^{0,q,1}(TM^0),\quad d"(\omega)=\xi_2(d_{01}\omega),
\end{eqnarray}
so we have 
\begin{eqnarray}
	d_{01}|_{\Omega^{0,q,0}(TM^0)}=d'+d".
\end{eqnarray}

\begin{proposition}
The operator $d'$ defined in relation (2.10) satisfies the relations:

a) $d'(\omega\wedge\theta)=d'\omega \wedge\theta +(-1)^q\omega\wedge d'\theta$, $(\forall)$$\omega\in \Omega^{0,q,0}(TM^0)$ and $\theta \in \Omega^{0,r,0}(TM^0)$.

b) $d'^2=0$.
\end{proposition}

\textit{Proof:}a) Let $\omega\in \Omega^{0,q,0}(TM^0)$ and $\theta \in \Omega^{0,r,0}(TM^0)$. It is known that $d_{01}(\omega\wedge\theta)=d_{01}\omega\wedge\theta +(-1)^q\omega\wedge d_{01}\theta$, and from relation (2.12), it follows 
	\[d'(\omega\wedge\theta)+d"(\omega\wedge\theta)=d'\omega\wedge\theta+d"\omega\wedge\theta +(-1)^q\omega\wedge d'\theta+(-1)^q\omega\wedge d"\theta.
\]
Considering the $(0,q+1,0)$-component in the both members, we have the desired result.

b) Let $\omega$ be a $(0,q,0)$-form. The definition (2.10) of the operator $d'$ says that $d'\omega=d_{01}\omega-\omega_0\wedge i_Zd_{01}\omega$. Hence we have
	\[d'^2\omega=d_{01}(d'\omega)-\omega_0\wedge i_Zd_{01}d'\omega=\omega_0\wedge d_{01}i_Zd_{01}\omega+\omega_0\wedge i_Z(d_{01}\omega\wedge i_Zd_{01}\omega-\omega_0\wedge d_{01}i_Zd_{01}\omega,
\]
where we used relations $d_{01}^2=0$, $d_{01}\omega_0=0$. Computing the last member in the above equalities, we obtain $d'^2=0$.

\begin{remark}
From the relation (2.11) ones can see that the operator $d"$ could not be composed to itself, but it can be proved the equality $d_{01}\circ d"+d"\circ d'=0$. 
\end{remark}

\begin{example}
a) For a $(0,1)$-form $\omega$, we have $\xi_1(\omega)=\omega-\omega(Z)\omega_0$, and $\xi_2(\omega)=\omega(Z)\omega_0$.

b) Let $f\in \Omega^0(TM^0)$, and $d_{01}f$ its foliated derivative, locally given by $d_{01}f=\frac{\partial f}{\partial y^i}\delta y^i$. Locally we have
	\[d"f=\xi_2(d_{01}f)=(d_{01}f)(Z)\omega_0=Z(f)\omega_0, 
\]
 
	\[d'f=d_{01}f-Z(f)\omega_0=\frac{\partial f}{\partial y^i}\delta y^i-y^i\frac{\partial f}{\partial y^i}\omega^0=\frac{\partial f}{\partial y^i}\theta_i,
\]
where $\theta_i$ are the $(0,1,0)$-forms defined in Example 2.1. Moreover, taking into account relation (1.10) and the fact $\sum_{i=1}^n t_i\theta_i=0$, it results that locally 
\begin{eqnarray}
	d'f=(X_if)\theta_i.
\end{eqnarray}
We have $d'y^j=(X_iy^j)\theta_i=\delta_i^j\theta_i-t_iZ(y^j)\theta_i=\theta_j-(t_i\theta_i)y^j=\theta_j$, so the $(0,1,0)$-forms $\theta_i$ are exactly the $d'$-derivatives of the local coordinates $y^i$, for all $i=\overline{1,n}$.

c) The $(0,2,0)$-forms $d'y^i\wedge d'y^j$ are $d'$-closed forms, for all $i,j=\overline{1,n}$.
\end{example}

Let us consider an arbitrar $(0,1)$-form on $TM^0$. It is locally given in $U$ by $\omega=\sum a_i\delta y^i$, with $a_i\in \Omega^0(U)$ such that in $U\cap \tilde{U}$,
	
\begin{eqnarray}
	\tilde{a}_{i_1}=\frac{\partial x^i}{\partial \tilde{x}^{i_1}}a_i.
\end{eqnarray}

From Proposition 2.2, $\omega$ is a $(0,1,0)$-form iff $i_Z\omega=0$, which is equivalent locally with $\sum a_i y^i=0$.
Then, locally we have
	\[\omega=\sum a_i \delta y^i=\sum a_i(d'y^i+y^i\omega_0)=\sum a_id'y^i+(\sum a_i y^i)\omega_0=\sum a_i d'y^i.
\]
  
Conversely, the expression locally given by $\sum a_id'y^i$, with functions $a_i$ satisfying (2.14) is a $(0,1,0)$-form because $d'y^i(Z)=0$, for all $i=\overline{1,n}$.

 \section{The d'-cohomology }
 
 \setcounter{equation}{0}
\begin{definition}
A function $f\in \Omega^0(TM^0)$ is called vertical Liouville basic if $d'f=0$. We denote by $\Sigma^0(TM^0)$ the space of such functions.
\end{definition}

The above definition and the Proposition 2.7 b) prove that the sequence
\begin{eqnarray}
	O\rightarrow \Sigma^0(TM^0)\stackrel{i}{\rightarrow} \Omega^0(TM^0)\stackrel{d'}{\rightarrow}\Omega^{0,1,0}(TM^0) \stackrel{d'}{\rightarrow}\Omega^{0,2,0}(TM^0)\stackrel{d'}{\rightarrow}...\stackrel{d'}{\rightarrow},
\end{eqnarray}
is a semiexact one.
\begin{definition}
We say that a $(0,q,0)$-form $\omega$ is $d'$-closed if $d'\omega=0$. We say that it is $d'$-exact if $\omega=d'\theta$ for some $\theta\in \Omega^{0,q-1,0}(TM^0)$. We denote by $Z^{0,q,0}(TM^0)$, $B^{0,q,0}(TM^0)$ the spaces of the $d'$-closed and $d'$-exacts $(0,q,0)$-forms, respectively.
\end{definition}
Taking into account that $d'^2=0$, we have the inclusion
	\[B^{0,q,0}(TM^0)\subset Z^{0,q,0}(TM^0).
\]
We call the \textit{d'-cohomology group of $TM^0$} the quotient group 
\begin{eqnarray}
	H_{d'}^{0,q,0}(TM^0)=\frac{Z^{0,q,0}(TM^0)}{B^{0,q,0}(TM^0)}.
\end{eqnarray}
This group is the de Rham group of the sequence (3.1).
\begin{theorem}
The operator $d'$ is satisfying a Poincare type lemma: Let $\omega \in \Omega^{0,q,0}(TM^0)$ be a d'-closed form. For every domain $U$ there is a $(0,q-1,0)$-form $\theta$ on $U$ such that locally $\omega=d'\theta$. 
\end{theorem}

\textit{Proof:} Let $\omega\in \Omega^{0,q,0}(TM^0)$ such that $d'\omega=0$. Then
	\[d_{01}\omega=d'\omega+d"\omega=d"\omega=\omega_0\wedge i_Zd_{01}\omega, 
\]so
	\[d_{01}\omega \equiv 0(mod \quad \omega_0).
\]Hence on the space $\omega_0=0$ we have $\omega$ $d_{01}$-exact. But the foliated derivative with respect to the vertical foliation satisfies a Poincare type lemma, so in every domain $U$ there is a vertical $(q-1)$-form $\theta$ such that
	\[\omega=d_{01}\theta +\lambda \wedge \omega_0, \quad \lambda\in \Omega^{0,q-1}(TM^0).
\]
 
Following proposition 2.3, $\theta=\omega_0\wedge i_Z\theta+\theta_1$, with $\theta_1=\xi_1(\theta)\in \Omega^{0,q-1,0}(TM^0)$. We obtain
	\[\omega=-\omega_0\wedge d_{01}i_Z\theta+d_{01}\theta_1+\lambda \wedge \omega_0.
\]
Here $\omega$ is a $(0,q,0)$-form, $\omega_0\wedge (d_{01}i_Z\theta+\lambda)$ is $(0,q-1,1)$-form and 
	\[d_{01}\theta_1=d'\theta_1+d"\theta_1\in \Omega^{0,q,0}(U)\oplus \Omega^{0,q-1,1}(U).
\]
 It results $\omega |_U=d'\theta_1$, q.e.d.

 Taking into account that $\Omega^0$ and the shaves of germs of $(0,q,0)$-forms are fins (see \cite{Mir1}, P.6.2, p.269), a consequence of the theorem 3.1 is the following:
 
\begin{proposition}
The sequence of shaves 
\begin{eqnarray}
	O\rightarrow \Sigma^0\stackrel{i}{\rightarrow} \Omega^0\stackrel{d'}{\rightarrow}\Omega^{0,1,0} \stackrel{d'}{\rightarrow}\Omega^{0,2,0}\stackrel{d'}{\rightarrow}...\stackrel{d'}{\rightarrow},
\end{eqnarray}
is a fine resolution of the sheaf $\Sigma^0$ of germs of vertical Liouville basic functions.
\end{proposition}
Now, by a well-known theorem of algebraic topology (see \cite{Mir1} T.3.6, p.205), we have the main result of this paper, a de Rham type theorem for the $d'$-cohomology:
\begin{theorem}
The q-dimensional Cech cohomology group of $TM^0$ with coefficients in $\Sigma^0$ is isomorphic with $H^{0,q,0}_{d'}(TM^0)$.
\end{theorem}

\textbf{Acknowledgment}: This work was supported by Contract with Sinoptix
No. 8441/2009.

\end{document}